\input{aipcheck}

\documentclass[
    ,final          % use final for the camera ready runs
  ]
  {aipproc}

\layoutstyle{8x11single}

\def\bR {\mathbf{R}}
\def\bS {\mathbf{S}}

\def\cC {\mathcal{C}}
\def\cD {\mathcal{D}}

\def\cL {\mathcal{L}}
\def\cM {\mathcal{M}}

\def\cQ {\mathcal{Q}}

\def\a {{\alpha}}

\def\de {{\delta}}

\def\eps {{\epsilon}}
\def\th {{\theta}}
\def\ka {{\kappa}}

\def\om {{\omega}}
\def\Om {{\Omega}}

\def\rstr {{\big |}}

\def\la {\langle}
\def\ra {\rangle}
\def \La {\bigg\langle}
\def \Ra {\bigg\rangle}

\def\Kn {\mathrm{Kn}}
\def\Sh {\mathrm{Sh}}
\def\Ma {\mathrm{Ma}}

\def\d {{\partial}}
\def\dd {\,\mathrm{d}}
\def\grad {{\nabla}}
\def\Dlt {{\Delta}}
\def\Div{\mathrm{div}}
\def\Tr{\mathrm{Tr}}

\newcommand{\ba}{\begin{aligned}}
\newcommand{\ea}{\end{aligned}}

\def\iint {\int\!\!\!\int}
\def\iiint {\int\!\!\!\int\!\!\!\int}
\def\tfrac#1#2 {{\textstyle\frac{#1}#2}}

\begin{document}

\title{From the Kinetic Theory of Gases to Continuum Mechanics}


\author{Fran\c cois Golse}{
  address={Ecole Polytechnique, Centre de Math\'ematiques Laurent Schwartz, 91128 Palaiseau Cedex, France}
}

\begin{abstract}
Recent results on the fluid dynamic limits of the Boltzmann equation based on the DiPerna-Lions theory of renormalized solutions are reviewed
in this paper, with an emphasis on regimes where the velocity field behaves to leading order like that of an incompressible fluid with constant 
density.
\end{abstract}

\maketitle

\rightline{\it In memory of Carlo Cercignani (1939--2010)}

\bigskip
%%%%%%%%%%%%%%%%%%%%%%%%%%%%%%%%%%%%%%%%%%%%
%% MAINMATTER
%%%%%%%%%%%%%%%%%%%%%%%%%%%%%%%%%%%%%%%%%%%%

Relating the kinetic theory of gases to their description by the equations of continuum mechanics is a question that finds its origins in the work of
Maxwell \citep{Maxwell66}. It was subsequently formulated by Hilbert as a mathematical problem --- specifically, as an example of his 6th problem
on the axiomatization of physics \citep{Hilbert00}. In Hilbert's own words ``Boltzmann's work on the principles of mechanics suggests the problem of 
developing mathematically the limiting processes which lead from the atomistic view to the laws of motion of continua''. Hilbert himself studied this
problem in \citep{Hilbert12} as an application of his theory of integral equations. The present paper reviews recent progress on this problem in the
past 10 years as a consequence of the DiPerna-Lions global existence and stability theory \citep{diPernaLions90} for solutions of the Boltzmann 
equation. This Harold Grad lecture is dedicated to the memory of Carlo Cercignani, who gave the first Harold Grad lecture in the 17th Rarefied 
Gas Dynamics Symposium, in Aachen (1990), in recognition of his outstanding influence on the mathematical analysis of the Boltzmann equation
in the past 40 years.

%%%%%%%%%%%%%%%%%%%%%%%%%%%%%%%%%%%%%%%%%%%%%%%%%%%%%%%%%%%%%%%%%%%%%%%%%
\section{The Boltzmann equation: formal structure}
%%%%%%%%%%%%%%%%%%%%%%%%%%%%%%%%%%%%%%%%%%%%%%%%%%%%%%%%%%%%%%%%%%%%%%%%%

In kinetic theory, the state of a monatomic gas at time $t$ and position $x$ is described by its velocity distribution function $F\equiv F(t,x,v)\ge 0$.
It satisfies the Boltzmann equation
$$
\d_tF+v\cdot\grad_xF=\cC(F)\,,
$$
where $\cC(F)(t,x,v):=\cC(F(t,x,\cdot))(v)$ is the Boltzmann collision integral defined for each continuous, rapidly decaying function $f\equiv f(v)$ 
by 
$$
\cC(f)(v):=\iint_{\bR^3\times\bS^2}(f(v')f(v'_*)-f(v)f(v_*))\tfrac{d^2}{2} |(v-v_*)\cdot\om|\dd v_*\dd \om\,,
$$
assuming that gas molecules behave as perfectly elastic hard spheres of diameter $d$. In this formula, we have denoted
\begin{equation}\label{CollPairFla}
v'\equiv v'(v,v_*,\om):=v\,-(v-v_*)\cdot\om\om\,,\qquad v'_*\equiv v'_*(v,v_*,\om):=v_*\!+(v-v_*)\cdot\om\om\,,\qquad|\om|=1\,.
\end{equation}
Molecular interactions more general than hard sphere collisions can be considered by replacing $\frac{d^2}{2}|(v-v_*)\cdot\om|$ with appropriate
collision kernels of the form $b(|v-v_*|,|\frac{v-v_*}{|v-v_*|}\cdot\om|)$. In this paper, we restrict our attention to the case of hard sphere collisions 
to avoid dealing with more technical conditions on the collision kernel.

%%%%%%%%%%%%%%%%%%%%%%%%%%%%%%%%%%%%%%%%%%%%%%%%%%%%%%%%%%%%%%%%%%%%%%%%%
\subsection{Properties of the collision integral}
%%%%%%%%%%%%%%%%%%%%%%%%%%%%%%%%%%%%%%%%%%%%%%%%%%%%%%%%%%%%%%%%%%%%%%%%%

While the collision integral is a fairly intricate mathematical expression, the formulas (\ref{CollPairFla}) entail remarkable symmetry properties. 
As a result, the collision integral satisfies, for each continuous, rapidly decaying $f\equiv f(v)$, the identities
\begin{equation}\label{CollCons}
\int_{\bR^3}\cC(f)\dd v=0\,,\quad\int_{\bR^3}\cC(f)v_k\dd v=0\,,\quad k=1,2,3,\quad\hbox{ and }\quad \int_{\bR^3}\cC(f)|v|^2\dd v=0\,.
\end{equation}
The first relation expresses the conservation of mass (or equivalently, of the number of particles) by the collision process, while the second and the 
third express the conservation of momentum and energy respectively. 

Perhaps the most important result on the collision integral is \textit{Boltzmann's H Theorem:} for  each continuous, rapidly decaying $f\equiv f(v)>0$
such that $\ln f$ has polynomial growth as $|v|\to+\infty$,
\begin{equation}\label{HThm}
\int_{\bR^3}\cC(f)\ln f\dd v\le 0\,,\qquad\hbox{ and }\quad\int_{\bR^3}\cC(f)\ln f\dd v=0\Leftrightarrow \cC(f)=0\Leftrightarrow f\hbox{ is a Maxwellian,}
\end{equation}
i.e. there exists $\rho,\th>0$ and $u\in\bR^3$ such that
\begin{equation}\label{Mawx}
f(v)=\cM_{(\rho,u,\th)}(v):=\frac{\rho}{(2\pi\th)^{3/2}}\exp\left(-\frac{|v-u|^2}{2\th}\right)\,.
\end{equation}

%%%%%%%%%%%%%%%%%%%%%%%%%%%%%%%%%%%%%%%%%%%%%%%%%%%%%%%%%%%%%%%%%%%%%%%%%
\subsection{Dimensionless variables}
%%%%%%%%%%%%%%%%%%%%%%%%%%%%%%%%%%%%%%%%%%%%%%%%%%%%%%%%%%%%%%%%%%%%%%%%%

Fluid dynamic limits are obtained as properties of solutions of the Boltzmann equation under appropriate scaling assumptions. We therefore
recast the Boltzmann equation in dimensionless variables, so as to identify the dimensionless parameters that control the scalings of the time
and space variables, following \citep{BGL2,SoneBook1}.

First we choose a macroscopic length scale $L$ (for instance the size of the container where the gas is enclosed, or of an object immersed in
the fluid, or the typical length scale on which the variation of macroscopic fluid quantities is observed), as well as a macroscopic observation 
time scale $T_o$ (i.e. the time scale on which the evolution of the fluid quantities is observed.)

We next define reference scales of density $\overline\rho$ and temperature $\overline\theta$ by setting
$$
\iint F\dd x\dd v=\overline\rho L^3\,,\qquad\iint vF \dd x\dd v=0\,,\qquad\iint\tfrac12 |v|^2F\dd x\dd v=\tfrac32 \overline\rho\overline\th\,.
$$
The collision time scale $T_c$ is defined in turn by
$$
\tfrac{d^2}{2} \iiint\cM_{(\overline{\rho},0,\overline{\th})}(v)\cM_{(\overline{\rho},0,\overline{\th})}(v_*)|(v-v_*)\cdot\om|\dd v\dd v_*\dd\om
	=\frac{\overline{\rho}}{T_c}\,,
$$
while the acoustic time scale is defined by $T_a=L/\sqrt{\overline{\th}}$. The dimensionless variables are $\hat t=t/T_o$, $\hat x=x/L$, and
$\hat v=v/{\sqrt{\overline\th}}$, while the dimensionless distribution function is $\hat F={\overline{\th}^{3/2}}F/{\overline\rho}$.

Introducing two dimensionless parameters, the Strouhal number $\Sh=T_a/T_o$ and the Knudsen number $\Kn=CT_c/T_a=2/d^2L\overline{\rho}$ 
with
$$
C=\tfrac1{{4\pi^2}} \iint|v-v_*|e^{-(|v|^2+|v_*|^2)/2}\dd v\dd v_*\,,
$$
we see that the Boltzmann equation in dimensionless variables takes the form
\begin{equation}\label{BoltzEq}
\Sh\,\d_{\hat t}\hat F+\hat v\cdot\grad_{\hat x}\hat F=\frac1{\Kn}\hat\cC(\hat F)\,,
\end{equation}
where the dimensionless collision integral is 
\begin{equation}\label{CollInt}
\hat\cC(\hat F)(\hat t,\hat x,\hat v)=\iint_{\bR^3\times\bS^2}(\hat F(\hat t,\hat x,\hat v')\hat F(\hat t,\hat x,\hat v'_*)
	-\hat F(\hat t,\hat x,\hat v)\hat F(\hat t,\hat x,\hat v_*))|(\hat v-\hat v_*)\cdot\om|\dd\hat v_*\dd\om\,.
\end{equation}
Obviously, the dimensionless collision integral $\hat\cC(\hat F)$ satisfies exactly the same properties as the original expression $\cC(F)$, i.e.
the conservation laws of mass, momentum and energy (\ref{CollCons}) and Boltzmann's H Theorem (\ref{HThm}).

Henceforth, we always consider the Boltzmann equation (\ref{BoltzEq}) in dimensionless variables, dropping all hats for notational simplicity. 
Thus, the conservation properties of the collision operator (\ref{CollCons}) imply that rapidly decaying (in $v$) solutions of the dimensionless 
Boltzmann equation (\ref{BoltzEq}) satisfy the following local conservation laws:
\begin{equation}\label{ConsLaws}
\left\{
\begin{array}{rl}
\Sh\,\displaystyle\d_t\int_{\bR^3}F\dd v+\Div_x\int_{\bR^3}vF\dd v=0&\quad\hbox{ (conservation of mass),}
\\	\\
\Sh\,\displaystyle\d_t\int_{\bR^3}vF\dd v+\Div_x\int_{\bR^3}v\otimes vF\dd v=0&\quad\hbox{ (conservation of momentum),}
\\	\\
\Sh\,\displaystyle\d_t\int_{\bR^3}\tfrac12 |v|^2F\dd v+\Div_x\int_{\bR^3}v\tfrac12 |v|^2F\dd v=0&\quad\hbox{ (conservation of energy).}
\end{array}
\right.
\end{equation}
Likewise, Boltzmann's H Theorem implies that solutions $F>0$ of the Boltzmann equation that are rapidly decaying while $\ln F$ has polynomial 
growth as $|v|\to+\infty$ satisfy the differential inequality
\begin{equation}\label{HDiffIneq}
\Sh\,\d_t\int_{\bR^3}F\ln F\dd v+\Div_x\int_{\bR^3}vF\ln F\dd v=\frac1{\Kn} \int_{\bR^3}\cC(f)\ln f\dd v\le 0\,.
\end{equation}

%%%%%%%%%%%%%%%%%%%%%%%%%%%%%%%%%%%%%%%%%%%%%%%%%%%%%%%%%%%%%%%%%%%%%%%%%
\section{The compressible Euler limit and Hilbert's expansion}
%%%%%%%%%%%%%%%%%%%%%%%%%%%%%%%%%%%%%%%%%%%%%%%%%%%%%%%%%%%%%%%%%%%%%%%%%

Whenever a gas evolves in a fluid dynamic regime (at the length scale $L$), the average time between successive collisions involving a typical
gas molecule is much smaller than the time necessary for an acoustic wave to travel a distance $L$. In other words, fluid dynamic regimes are
characterized by the condition $T_c\ll T_a$, or equivalently by he condition $\Kn\ll 1$.

In \citep{Hilbert12}, Hilbert studied the Boltzmann equation (\ref{BoltzEq}) in the asymptotic regime defined by $\Kn=\eps\ll 1$ and $\Sh=1$. His
idea was to seek the solution $F_\eps$ of
\begin{equation}
\label{BoltzEulScal}
\d_tF_\eps+v\cdot\grad_xF_\eps=\frac1\eps\cC(F_\eps)
\end{equation}
as a formal power series in $\eps$ with smooth coefficients --- known as \textit{Hilbert's expansion}:
\begin{equation}
\label{HilbSer}
F_\eps(t,x,v)=\sum_{n\ge 0}\eps^nf_n(t,x,v)\,,\qquad\hbox{ with }f_n\hbox{ smooth in }t,x,v\,,\,\,\hbox{ for each }n\ge 0\,.
\end{equation}
He found that the leading order term in that expansion is of the form 
$$
f_0(t,x,v)=\cM_{(\rho,u,\th)(t,x)}(v)\,,
$$
where $(\rho,u,\th)$ is a solution of the compressible Euler system
\begin{equation}\label{ComprEul}
\left\{
\begin{array}{r}
\d_t\rho+\Div_x(\rho u)=0\,,
\\	\\
\rho(\d_tu+u\cdot\grad_xu)+\grad_x(\rho\th)=0\,,
\\	\\
\d_t\th+u\cdot\grad_x\th+\tfrac23\th\Div_xu=0\,.
\end{array}
\right.
\end{equation}
Caflisch \citep{Caflisch80} succeeded in turning Hilbert's formal result into a rigorous statement bearing on solutions of the Boltzmann equation,
by using a truncated variant of the Hilbert expansion above. Specifically, given a smooth solution $(\rho,u,\th)$ of the compressible Euler system
on some finite time interval $[0,T)$, he constructs a family of solutions of the Boltzmann equation that converges to $\cM_{(\rho,u,\th)}$ uniformly
in $t\in[0,T)$ as $\eps\to 0$. Before Caflisch's result, Nishida \citep{Nishida78} had proposed another proof of the compressible Euler limit of the
Boltzmann equation under more stringent regularity assumptions, viz. analyticity, using some abstract variant of the Cauchy-Kowalewski theorem.

One striking advantage of the Hilbert expansion is its versatility, abundantly illustrated by the great diversity of physically meaningful applications
to be found in the work of Sone \citep{SoneBook1,SoneBook2}. However, there are some serious difficulties with the Hilbert expansion, some of 
which can be treated with adequate mathematical techniques. First, the radius of convergence of the Hilbert power series is $0$ in general, so that 
essentially all mathematical arguments based on Hilbert's expansion use a truncated variant thereof. In general, truncated Hilbert expansions are 
not everywhere nonnegative, and are not exact solutions of the Boltzmann equation. One obtains exact solutions of the Boltzmann equation by 
adding to the truncated Hilbert expansion some appropriate remainder term, satisfying a variant of the Boltzmann equation that becomes weakly 
nonlinear for small enough $\eps$ (see for instance \citep{Caflisch80,dMEL,ArkNouCouette}.) The truncated Hilbert expansion with the remainder 
term so constructed is a rigorous, pointwise asymptotic expansion (meaning that $\eps^{-n}|F_\eps-(f_0+\eps f_1+\ldots+\eps^nf_n)|\to 0$ 
pointwise in $(t,x,v)$) of the solution $F_\eps$ of (\ref{BoltzEulScal}) as $\eps\to 0$. Another difficulty in working with Hilbert's expansion, even 
truncated at some finite order, is that $f_n=O(|\grad^n_{\hat t,\hat x}f_0|)$ for each $n\ge 0$. Since generic solutions of the compressible Euler 
system lose regularity in finite time \citep{Sideris87}, truncated Hilbert expansions make sense on finite time intervals only. For instance, if a 
solution $(\rho,u,\th)$ of the compressible Euler system involves a shock wave, only the $0$-th order term in the associated Hilbert expansion, 
i.e. $f_0(t,x,v)=\cM_{(\rho,u,\th)(t,x)}(v)$ is well defined. In general, if the geometric structure and the position of the singularities in the solution 
of the hydrodynamic equations are known precisely, one can bypass this difficulty by adding to the truncated Hilbert expansion appropriate 
boundary layer terms. If the structure of these singularities is unknown, or one does not even know whether the hydrodynamic solution is smooth,
one cannot use the Hilbert expansion.

%%%%%%%%%%%%%%%%%%%%%%%%%%%%%%%%%%%%%%%%%%%%%%%%%%%%%%%%%%%%%%%%%%%%%%%%%
\section{Global existence theory for the Boltzmann equation}
%%%%%%%%%%%%%%%%%%%%%%%%%%%%%%%%%%%%%%%%%%%%%%%%%%%%%%%%%%%%%%%%%%%%%%%%%

To avoid the various shortcomings of the Hilbert expansion method, one needs a theory of global solutions for the Boltzmann equation based
on the only estimates that are uniform in $\Kn$ as $\Kn\to 0$. These estimates are those deduced from the conservation laws (\ref{CollCons})
and Boltzmann's H Theorem (\ref{HThm}), or from their differential formulations (\ref{ConsLaws})-(\ref{HDiffIneq}).

Henceforth, we are concerned with solutions of the Boltzmann equation for a gas filling the Euclidian space $\bR^3$ and at equilibrium at infinity. 
By Galilean invariance and with a convenient choice of units, we can assume without loss of generality, that this equilibrium state at infinity is the Maxwellian $\cM_{(1,0,1)}$, denoted by $M$ in the sequel. In other words, we seek the solution of
\begin{equation}\label{BoltzEqM}
\left\{
\begin{array}{l}
\displaystyle\Sh\,\d_tF+v\cdot\grad_xF=\frac1{\Kn}\cC(F)\,,\quad (x,v)\in\bR^3\times\bR^3\,,\,\,t>0\,,
\\	\\
F(t,x,v)\to M\quad\hbox{as }|x|\to+\infty\,,
\\	\\
F\rstr_{t=0}=F^{in}\,.
\end{array}
\right.
\end{equation}

A convenient quantity measuring the distance between two distribution functions in the context of the Boltzmann equation is the \textit{relative
entropy}: for $F\equiv F(x,v)\ge 0$ and $G\equiv G(x,v)>0$ a.e. in $(x,v)\in\bR^3\times\bR^3$,
\begin{equation}\label{RelEntr}
H(F|G):=\iint_{\bR^3\times\bR^3}(F\ln(F/G)-F+G)(x,v)\dd x\dd v\,.
\end{equation}
Notice that $a\ln(a/b)-a+b\ge 0$ for each $a\ge 0$ and $b>0$, with equality if and only if $a=b$. Hence the integrand is a nonnegative measurable
function and $H(F|G)=0$ if and only if $F=G$ a.e. on $\bR^3\times\bR^3$. 

Since $\ln M=-\frac32\ln(2\pi)-\frac12|v|^2$, a formal argument based on the local conservation laws (\ref{ConsLaws}) and the differential inequality 
(\ref{HDiffIneq}) shows that any classical solution $F$ of (\ref{BoltzEqM}) with appropriate decay as $|v|\to+\infty$ satisfies
$$
\Sh\,\d_t\int_{\bR^3}(F\ln(F/M)-F+M)\dd v+\Div_x\int_{\bR^3}v(F\ln(F/M)-F+M)\dd v\le 0\,.
$$
Integrating in $x$ both sides of this inequality and assuming that $F\to M$ fast enough as $|x|\to+\infty$, we conclude that 
\begin{equation}\label{EntrBnd}
\sup_{t\ge 0}H(F(t,\cdot,\cdot)|M)\le H(F^{in}|M)\,,\quad\hbox{ and }\quad
	\int_0^\infty\iint_{\bR^3\times\bR^3}-\cC(F)\ln F\dd s\dd x\dd v\le\Sh\,H(F^{in}|M)\,.
\end{equation}

Observe that the collision integral $\cC(F)$ acts as a nonlocal integral operator analogous to a convolution in the $v$ variable and as a pointwise
product in the $x$ variable. The fact that $\cC(F)$ is quadratic in $F$ while $H(F|M)$ is ``essentially homogeneous of degree $1$ as $F\gg 1$''
suggests that $\cC(F)$ may not be defined for all nonnegative measurable functions $F$ satisfying the entropy bound (\ref{EntrBnd}) above. Yet, 
for each measurable $F>0$ on $\bR^3\times\bR^3$ , one has 
$$
\iint_{|x|+|v|\le r}\frac{|\cC(F)|}{\sqrt{1+F}}\dd x\dd v\le C\iint_{|x|\le r}(-\cC(F)\ln F+(1+|v|^2)F)\dd x\dd v
$$
so that $\cC(F)/\sqrt{1+F}\in L^1_{loc}(\bR_+\times\bR^3\times\bR^3)$, i.e. is locally integrable in $(t,x,v)$. This suggests dividing both sides of the 
Boltzmann equation by $\sqrt{1+F}$, thereby leading to the notion of \textit{renormalized solution}.

\smallskip
\noindent
\textbf{Definition.} (DiPerna-Lions \citep{diPernaLions90})
A renormalized solution relative to $M$ of the Boltzmann equation is a nonnegative function $F\in C(\bR_+,L^1_{loc}(\bR^3\times\bR^3))$ 
satisfying $H(F(t)|M)<+\infty$ for each $t\ge 0$ and
$$
M(\Sh\,\d_t+v\cdot\grad_x)\Gamma(F/M)=\frac1{\Kn}\cC(F)\Gamma'(F/M)
$$
in the sense of distributions on $\bR_+^*\times\bR^3\times\bR^3$, for each $\Gamma\in C^1(\bR_+)$ satisfying $\Gamma'(Z)\le C/\sqrt{1+Z}$.

\smallskip
With this notion of solution, one can prove the global existence and weak stability of solutions of the Cauchy problem for the Boltzmann equation,
with initial data that are not necessarily small perturbations of either the vacuum state or of a Maxwellian equilibrium.

\smallskip
\noindent
\textbf{Theorem.} (DiPerna-Lions-Masmoudi \citep{diPernaLions90,Lions94,LionsMasmoudi2})
For each measurable $F^{in}\ge 0$ a.e. on $\bR^3\times\bR^3$ satisfying the condition $H(F^{in}|M)<+\infty$, there exists a renormalized solution 
of the Boltzmann equation (\ref{BoltzEqM}) with initial data $F^{in}$. This solution $F$ satisfies 
\begin{equation}\label{ApproxConsLaws1}
\left\{
\begin{array}{r}
\displaystyle\Sh\,\d_t\int_{\bR^3}F\dd v+\Div_x\int_{\bR^3}vF\dd v=0\,,
\\
\displaystyle\Sh\,\d_t\int_{\bR^3}vF\dd v+\Div_x\int_{\bR^3}v\otimes vF\dd v+\Div_xm=0\,,
\end{array}
\right.
\end{equation}
where $m=m^T\ge 0$ is a matrix-valued Radon measure, and the entropy inequality
\begin{equation}\label{EntrIneq}
\Sh\,H(F(t,\cdot,\cdot)|M)+\Sh\int_{\bR^3}\Tr(m(t))-\int_0^t\!\!\iint_{\bR^3\times\bR^3}\cC(F)\ln F(s,x,v)\dd s\dd x\dd v\le\Sh\,H(F^{in}|M)\,,\quad t>0\,.
\end{equation}

\smallskip
A classical solution of the Boltzmann equation with appropriate decay as $|v|\to+\infty$ would satisfy all these properties with $m=0$; besides the 
entropy inequality is a weakened variant of Boltzmann's H Theorem --- which would imply that this inequality is in fact an equality. 

\smallskip
The main advantage of the notion of renormalized solutions is that a) such solutions always exist for each initial data with finite relative entropy 
with respect to $M$, and b) such solutions are weakly stable, in the sense that if a sequence $(F_n)_{n\ge 0}$ of renormalized solutions of the
Boltzmann equation converges to $F$ in the sense of distributions and satisfies $H(F_n\rstr_{t=0}|M)\le C$ for all $n\ge 0$, where $C$ is some
positive constant, then $F$ is also a renormalized solution of the Boltzmann equation, satisfying (\ref{ApproxConsLaws1}) and (\ref{EntrIneq}). 
Unfortunately, there is no uniqueness theorem for this notion of solution, so that a renormalized solution of the Boltzmann equation is not 
completely determined by its initial data. But if the Cauchy problem for the Boltzmann equation has a classical solution $F$, each renormalized 
solution of the Boltzmann equation with the same initial data as $F$ coincides with $F$ a.e. in $(t,x,v)$ (see \citep{LionsKyoto93}.)

%%%%%%%%%%%%%%%%%%%%%%%%%%%%%%%%%%%%%%%%%%%%%%%%%%%%%%%%%%%%%%%%%%%%%%%%%
\section{Fluid dynamic limits of the Boltzmann equation}
%%%%%%%%%%%%%%%%%%%%%%%%%%%%%%%%%%%%%%%%%%%%%%%%%%%%%%%%%%%%%%%%%%%%%%%%%

As explained above, all fluid dynamic limits of the Boltzmann equation are characterized by the scaling condition $\Kn\ll 1$: hence we set 
$\Kn=\eps$ throughout the present section.

Besides, all the fluid dynamic limits considered in this paper correspond with weakly nonlinear regimes at the kinetic level --- which does not 
imply that the nonlinearities are weak at the macroscopic level. Such regimes have been systematically explored by Sone at the formal level 
(see \citep{SoneBook1} and the references therein), by using the Hilbert expansion method. In other words, the distribution functions $F$ 
considered are small perturbations of the Maxwellian state $M$ at infinity. Henceforth, we denote by $\de_\eps\ll 1$ the order of magnitude 
of the difference $F-M$. A typical example of such a distribution function is $F(t,x,v)=\cM_{(1,\de_\eps u(t,x),1)}(v)$, since 
$\cM_{(1,\de_\eps u(t,x),1)}(v)=M(v)(1+\de_\eps u(t,x)\cdot v+O(\de_\eps^2))$. 
In this example, the distribution function $F$ defines a velocity field $u_F$ and a temperature field $\th_F$ by the formulas
$$
u_F:=\frac{\displaystyle\int_{\bR^3}vF\dd v}{\displaystyle\int_{\bR^3}F\dd v}=\de_\eps u
	\qquad\hbox{ and }\quad\th_F:=\frac{\displaystyle\int_{\bR^3}|v-u_F|^2F\dd v}{3\displaystyle\int_{\bR^3}F\dd v}=1\,.
$$
Introducing the speed of sound $c_F:=\sqrt{5\th_F/3}$, we see that the Mach number $\Ma:=u_F/c_F=\de_\eps u$, so that the scaling parameter
$\de_\eps$ can be thought of as the (order of magnitude of the) Mach number.

%%%%%%%%%%%%%%%%%%%%%%%%%%%%%%%%%%%%%%%%%%%%%%%%%%%%%%%%%%%%%%%%%%%%%%%%%
\subsection{The acoustic limit}
%%%%%%%%%%%%%%%%%%%%%%%%%%%%%%%%%%%%%%%%%%%%%%%%%%%%%%%%%%%%%%%%%%%%%%%%%

The acoustic limit is the linearized variant of the compressible Euler limit considered by Hilbert himself. 

\smallskip
\noindent
\textbf{Theorem.} (Golse, Jiang, Levermore, Masmoudi \citep{GolseLvrmr2002, JingLvrmrMasmoudi})
Let $\Kn=\eps$, $\Ma=\de_\eps=O(\sqrt{\eps})$ and $\Sh=1$. For each $\rho^{in},u^{in},\th^{in}\in L^2(\bR^3)$, let $F_\eps$ be a family of 
renormalized solutions of the Boltzmann equation (\ref{BoltzEqM}) with initial data
$$
F^{in}_\eps=\cM_{(1+\de_\eps\rho^{in},\de_\eps u^{in},1+\de_\eps\th^{in})}\,.
$$
Then, in the limit as $\eps\to 0$, 
$$
\frac1{\de_\eps}\int_{\bR^3}(F_\eps(t,x,v)-M(v))(1,v,\tfrac13|v|^2-1)\dd v\to(\rho,u,\th)(t,x)
$$
in $L^1_{loc}(\bR_+^*\times\bR^3)$, where $(\rho,u,\th)$ is the solution of the acoustic system
$$
\left\{
\begin{array}{ll}
\d_t\rho+\Div_xu=0\,,&\quad\rho\rstr_{t=0}=\rho^{in}\,,
\\
\d_tu+\grad_x(\rho+\th)=0\,,&\quad u\rstr_{t=0}=u^{in}\,,
\\
\d_t\th+\tfrac23 \Div_xu=0\,,&\quad\th\rstr_{t=0}=\th^{in}\,.
\end{array}
\right.
$$

While the result in \citep{GolseLvrmr2002} holds for the most general class of molecular interactions satisfying some angular cutoff assumption
in the sense of Grad \citep{Grad63} (in fact, a much weaker version of Grad's assumption \citep{JingLvrmrMasmoudi,LvrmrMasmoudi2010}), 
an earlier contribution of the same authors with Bardos \citep{BGL3}  introduced a key new idea in the derivation of hydrodynamic limits of 
renormalized solutions of the Boltzmann equation and treated the case of bounded collision kernels (e.g. cutoff Maxwell molecules).

%%%%%%%%%%%%%%%%%%%%%%%%%%%%%%%%%%%%%%%%%%%%%%%%%%%%%%%%%%%%%%%%%%%%%%%%%
\subsection{The incompressible Euler limit}
%%%%%%%%%%%%%%%%%%%%%%%%%%%%%%%%%%%%%%%%%%%%%%%%%%%%%%%%%%%%%%%%%%%%%%%%%

It is a well-known fact that, in the low Mach number limit, the flow of an inviscid fluid can be approximately decomposed into its acoustic and
vortical modes, whose interaction vanishes with the Mach number. The result below explores the counterpart for vortical modes of the acoustic 
limit of the Boltzmann equation. Because of the low Mach number scaling, vortical modes evolve on a longer time scale than acoustic modes, consistently with the fact that the conditions $\grad_x(\rho+\th)=0$ and $\Div_xu=0$ characterize the equilibrium points of the acoustic system.

\smallskip
\noindent
\textbf{Theorem.} (Saint-Raymond \citep{SR2003})
Let $\Kn=\eps$, and $\Sh=\Ma=\de_\eps=\eps^\a$ with $0<\a<1$. Let $u^{in}\in H^3(\bR^3)$\footnote{The notation $H^m(\bR^n)$
designates the Sobolev space of square integrable functions on $\bR^n$ whose partial derivatives of order $\le m$ in the sense of distributions
are square integrable functions on $\bR^n$. A vector field is said to belong to $H^m(\bR^n)$ if all its components belong to $H^m(\bR^n)$.}
satisfy $\Div u^{in}=0$, and let $u\in C([0,T];H^3(\bR^3))$ be a solution of the incompressible Euler equations
$$
\left\{
\begin{array}{l}
\d_tu+u\cdot\grad_xu+\grad_xp=0\,,\quad\Div_xu=0\,,
\\	
u\rstr_{t=0}=u^{in}\,.
\end{array}
\right.
$$
Let $F_\eps$ be a family of renormalized solutions of the Boltzmann equation (\ref{BoltzEqM}) with initial data
$$
F^{in}_\eps=\cM_{(1,\de_\eps u^{in},1)}\,.
$$
Then, in the limit as $\eps\to 0$, one has
$$
\frac1{\de_\eps}\int_{\bR^3}vF_\eps(t,x,v)\dd v\to u(t,x)\quad\hbox{ in }L^\infty([0,T];L^1_{loc}(\bR^3))\,.
$$

The proof of this result is based on the relative entropy method, described in the next section. Actually, there had been precursors of this theorem 
due to the author \citep{BouGolPul2000} and to Lions-Masmoudi \citep{LionsMasmoudi2}, where the relative entropy method was introduced for this
type of problem. Unfortunately, the statements in \citep{BouGolPul2000,LionsMasmoudi2} rested on extra assumptions on the family of solutions 
of the Boltzmann equation that remain unverified. These assumptions were removed by some clever argument in \citep{SR2003}, which therefore 
contains the first complete proof of the theorem above. 

%%%%%%%%%%%%%%%%%%%%%%%%%%%%%%%%%%%%%%%%%%%%%%%%%%%%%%%%%%%%%%%%%%%%%%%%%
\subsection{The Stokes limit}
%%%%%%%%%%%%%%%%%%%%%%%%%%%%%%%%%%%%%%%%%%%%%%%%%%%%%%%%%%%%%%%%%%%%%%%%%

We continue our exploration of vortical modes with the Stokes limit of the Boltzmann equation. The scaling is weakly nonlinear at the macroscopic
level of description, and the time scale is chosen so as to keep track of entropy production in the fluid dynamic limit. 

\noindent
\smallskip
\textbf{Theorem.} (Golse, Levermore, Masmoudi \citep{GolseLvrmr2002,LvrmrMasmoudi2010})
Let $\Kn=\Sh=\eps$, and $\Ma=\de_\eps=o(\eps)$. For each $(u^{in},\th^{in})\in L^2\times L^\infty(\bR^3)$ such that $\Div_xu^{in}=0$ and each
$\eps\in(0,\|u^{in}\|_{L^\infty})$, let  $F_\eps$ be a family of renormalized solutions of the Boltzmann equation (\ref{BoltzEqM}) with initial data
$$
F^{in}_\eps=\cM_{(1-\de_\eps\th^{in},\de_\eps u^{in},1+\de_\eps\th^{in})}\,.
$$
Then, in the limit as $\eps\to 0$, one has
$$
\frac1{\de_\eps}\int_{\bR^3}(F_\eps(t,x,v)-M(v))(v,\tfrac13|v|^2-1)\dd v\to(u,\th)(t,x)\hbox{ in }L^1_{loc}(\bR_+\times\bR^3\times\bR^3)\,,
$$
where $(u,\th)$ is a solution of the Stokes-Fourier system 
$$
\left\{
\begin{array}{ll}
\d_tu+\grad_xp=\nu\Dlt_xu\,,\quad\Div_x u=0\,,&\qquad u\rstr_{t=0}=u^{in}\,,
\\
\d_t\th=\tfrac25 \ka\Dlt_x\th\,,&\qquad\th\rstr_{t=0}=\th^{in}\,.
\end{array}
\right.
$$

\smallskip
The viscosity and heat conductivity in this theorem are given by the formulas (equivalent to the usual ones  in \citep{SoneBook1}):
\begin{equation}\label{nu-ka}
\nu=\tfrac1{5} \cD^*(v\otimes v-\tfrac13 |v|^2I)\,,\qquad\ka=\tfrac23 \cD^*(\tfrac12 (|v|^2-5)v)\,,
\end{equation}
where $\cD^*$ denotes the Legendre dual of the Dirichlet form $\cD$ of the collision operator linearized about $M$, i.e.
$$
\cD(\Phi):=\tfrac12 \iiint_{\bR^3\times\bR^3\times\bS^2}|\Phi+\Phi_{*}-\Phi'-\Phi'_{*}|^2|(v-v_*)\cdot\om|MM_*\dd v\dd v_*\dd\om\,.
$$

\smallskip
The fluid dynamic model obtained in the statement above is the Stokes-Fourier system; notice that the motion and temperature equations are
decoupled in the absence of an external force field deriving from a potential. Previously Lions and Masmoudi \citep{LionsMasmoudi2} arrived
at the particular case of the statement above corresponding to an initial data for which $\th^{in}=0$, leading to the motion equation only, i.e. the
evolution Stokes equation. For want of a better control of the high speed tails of the distribution function, their argument cannot be generalized 
to obtain the Stokes-Fourier system presented above. The proof in \citep{GolseLvrmr2002} uses a different idea originating from \citep{BGL3}.

%%%%%%%%%%%%%%%%%%%%%%%%%%%%%%%%%%%%%%%%%%%%%%%%%%%%%%%%%%%%%%%%%%%%%%%%%
\subsection{The incompressible Navier-Stokes limit}
%%%%%%%%%%%%%%%%%%%%%%%%%%%%%%%%%%%%%%%%%%%%%%%%%%%%%%%%%%%%%%%%%%%%%%%%%

Finally, we remove the weakly nonlinear scaling assumption at the macroscopic level of description, while keeping entropy production effects 
at leading order, and obtain the incompressible Navier-Stokes equations as a fluid dynamic limit of the Boltzmann equation.

\smallskip
\noindent
\textbf{Theorem.} (Golse, Saint-Raymond \citep{GSR2004,GSR2009})
Let $\Kn=\Sh=\Ma=\de_\eps=\eps$. For each $(u^{in},\th^{in})\in L^2\times L^\infty(\bR^3)$ such that $\Div_xu^{in}=0$, let $F_\eps$ be a family 
of renormalized solutions of the Boltzmann equation (\ref{BoltzEqM}) with initial data
$$
F^{in}_\eps=\cM_{(1-\eps\th^{in},\eps u^{in},1+\eps\th^{in})}\,,
$$
for each $\eps\in(0,\|u^{in}\|_{L^\infty})$. There exists at least one subsequence $\eps_n\to 0$ such that
$$
\frac1{\eps_n}\int_{\bR^3}(\!F_{\eps_n}\!(t,x,v)\!-\!M(v))\!\left(v,\tfrac13 |v|^2\!-\!1\right)\dd v\to(u,\th)(t,x)
	\hbox{ in weak-}L^1_{loc}(\bR_+\times\bR^3\times\bR^3)\,,
$$
where $(u,\th)$ is a ``Leray solution'' of the Navier-Stokes-Fourier system with viscosity $\nu$ and heat conductivity $\ka$ given by formula 
(\ref{nu-ka}):
$$
\left\{
\begin{array}{ll}
\d_tu+\Div_x(u\otimes u)+\grad_xp=\nu\Dlt_xu\,,\quad\Div_xu=0\,,&\qquad u\rstr_{t=0}=u^{in}\,,
\\
\d_t\th+\Div_x(u\th)=\tfrac25 \ka\Dlt_x\th\,,&\qquad \th\rstr_{t=0}=\th^{in}\,.
\end{array}
\right.
$$

\smallskip
Let us briefly recall the notion of \textit{Leray solution} of the Navier-Stokes-Fourier system. In \citep{Leray1934} (arguably one of the most important
papers in the modern theory of partial differential equations), Leray defined a convenient notion of weak solution of the Navier-Stokes equations,
and proved that, in space dimension $3$, any initial velocity field with finite kinetic energy launches at least one such solution defined for all 
times. Leray solutions are not known to be uniquely defined by their initial data; however, if an initial data launches a smooth solution, all Leray 
solutions with the same initial data must coincide with that smooth solution. At the time of this writing, it is yet unknown (and a major open 
problem in the analysis of partial differential equations) whether Leray solutions launched by any smooth initial data remain smooth for all times. 
Thus, we do not know whether different subsequences $\eps_n\to 0$ in the theorem above lead to the same Leray solution $(u,\th)$ of the 
Navier-Stokes-Fourier system in general.

\smallskip
A Leray solution of the Navier-Stokes-Fourier system above is a pair $(u,\th)$ consisting of a velocity field $u$ and a temperature field $\th$, both
continuous on $\bR_+$ with values in $L^2(\bR^3)$ equipped with its weak topology, that solves the Navier-Stokes-Fourier system in the sense of 
distributions, satisfies the initial condition, and verifies the \textit{Leray inequality}:
\begin{equation}\label{LerayIneq}
\tfrac12 \int_{\bR^3}(|u|^2+\tfrac52 |\th|^2)(t,x)\dd x+\int_0^t\int_{\bR^3}(\nu|\grad_xu|^2+\ka|\grad_x\th|^2)(s,x)\dd x\dd s
	\le\tfrac12 \int_{\bR^3}(|u^{in}|^2+\tfrac52 |\th^{in}|^2)(t,x)\dd x\,.
\end{equation}
The Leray inequality is an equality for classical solutions of the Navier-Stokes equations, exactly as the DiPerna-Lions entropy inequality 
(\ref{EntrIneq}) is an equality for classical solutions of the Boltzmann equation. This indicates that the Leray existence theory for the Navier-Stokes 
equations and the DiPerna-Lions existence theory for the Boltzmann equation are parallel theories. The theorem above explains how these 
theories are related in the hydrodynamic limit.

Partial results on this theorem have been obtained by Lions-Masmoudi \citep{LionsMasmoudi1}. While the reference \citep{GSR2004} treated the 
case of bounded collision kernels, the theorem above was later extended to all hard cutoff potentials in the sense of Grad --- which includes the
case of hard spheres considered in this paper --- in \citep{GSR2009}. The arguments in \citep{GSR2004,GSR2009} have been recently refined
by Levermore and Masmoudi \citep{LvrmrMasmoudi2010} to encompass both soft as well as hard potentials, under a cutoff assumption more 
general than that proposed by Grad in \citep{Grad63}.

While these results bear on the most general case of renormalized solutions without restrictions on the size of initial data in space dimension
$3$, the Navier-Stokes limit of the Boltzmann equation had previously been obtained in the case of global smooth solutions for small initial 
data by Bardos and Ukai \citep{BardosUkai91}. The Navier-Stokes limit of the Boltzmann equation had also been established on finite time
intervals by adapting the Caflisch method based on Hilbert truncated expansions, by DeMasi, Esposito and Lebowitz \citep{dMEL}.

\medskip
The fluid dynamic limits discussed in this section can therefore be summarized as in table \ref{FluidDynLimTable}. Notice that these limits have
been established for molecular interactions more general than hard sphere collisions; see the references listed in the statements of the various
theorems above for the conditions on the collision kernel $b(v-v_*,\om)$. All these results assume some angular cutoff on the collision kernel as 
proposed by Grad \citep{Grad63} --- or slightly more general, as in \citep{LvrmrMasmoudi2010}. 

\smallskip
More importantly, some of the conditions bearing on the parameters $\Kn$, $\Ma$ and $\Sh$ may be not optimal. Formal arguments suggest that
the acoustic limit should hold whenever $\de_\eps\ll 1$ instead of $\de_\eps\ll\sqrt\eps$, while the incompressible Euler limit should hold under
the weaker condition $\de_\eps\gg\eps$ instead of $\de_\eps=\eps^\a$ with $0<\a<1$.

\begin{table}[!t]
\label{FluidDynLimTable}
\begin{tabular}{|c|c|c|}
\hline
%\multicolumn{4}{|c|}{}\\
\multicolumn{3}{|c|}{\bf Boltzmann equation $\Kn=\eps\ll 1$}\\
%\multicolumn{4}{|c|}{}\\
%\hline
%\multicolumn{4}{|c|}{}\\
%\multicolumn{3}{|c|}{\bf von Karman relation $\Ma/\Kn=\Rey$}\\
%\multicolumn{4}{|c|}{}\\
\hline
%& & \\
\hspace{.2cm} $\Ma$ \hspace{.2cm} & \hspace{.2cm} $\Sh$ \hspace{.2cm}& \hspace{.2cm} Fluid dynamic limit \hspace{.2cm} \\
%& & \\
\hline
%\hline
%& & \\
$\de_\eps\ll\sqrt\eps $ & $1$ & Acoustic system\\
%& & \\
\hline
%& & \\
$\de_\eps\ll\eps$ & $\eps$ & Stokes-Fourier system\\
%& & \\
\hline
%& & \\
$\de_\eps=\eps^\a\,,\,\,0<\a<1 $ & $\de_\eps$ & Incompressible Euler equations\\
%& & \\
\hline
%& & \\
$\eps$ & $\eps$ & Incompressible Navier-Stokes equations\\
%& & \\
\hline
\end{tabular}
\caption{Fluid dynamic limits of the Boltzmann equation, depending on the dimensionless parameters $\Kn$, $\Ma$ and $\Sh$.} 
\end{table}

Let us conclude this section with an important remark on the physical meaning of the ``incompressible'' fluid dynamic limits of the Boltzmann 
equation. What is proved in the last three theorems is that, to leading order, the velocity field $u$ satisfies the same equations  as the velocity 
field of an incompressible fluid with constant density. This does not mean that the gas is incompressible in that regime. Also, in the case of an 
incompressible fluid with the same heat capacity and heat conductivity as the gas, the diffusion term in the equation for the temperature field 
would be multiplied by $5/3$. This difference comes from the work of the pressure: see the discussion in footnotes 6 on p. 93 in \citep{SoneBook1} 
and 43 on p. 107 of \citep{SoneBook2}, together with section 3.7.2 in \citep{SoneBook2}. 

Likewise, the inequality (\ref{LerayIneq}) was written by Leray in \citep{Leray1934} with $\th\equiv 0$. For an incompressible fluid with constant 
density $\overline\rho$, the quantity $\frac12\int\overline\rho|u(t,x)|^2\dd x$ is the kinetic energy of the fluid at time $t$, and the Leray inequality is 
interpreted as a statement on the dissipation of energy in the fluid. The meaning of (\ref{LerayIneq}) with $\th\not\equiv0$ is obviously different, 
since the quantity $\frac12\int(|u(t,x)|^2+\frac52\th(t,x)^2)\dd x$ is not the total energy of the gas at time $t$.

%%%%%%%%%%%%%%%%%%%%%%%%%%%%%%%%%%%%%%%%%%%%%%%%%%%%%%%%%%%%%%%%%%%%%%%%%
\section{Mathematical tools for the hydrodynamic limit}
%%%%%%%%%%%%%%%%%%%%%%%%%%%%%%%%%%%%%%%%%%%%%%%%%%%%%%%%%%%%%%%%%%%%%%%%%

%%%%%%%%%%%%%%%%%%%%%%%%%%%%%%%%%%%%%%%%%%%%%%%%%%%%%%%%%%%%%%%%%%%%%%%%%
\subsection{The linearized collision integral}
%%%%%%%%%%%%%%%%%%%%%%%%%%%%%%%%%%%%%%%%%%%%%%%%%%%%%%%%%%%%%%%%%%%%%%%%%

In all the fluid dynamic limits considered in the previous section, the solution $F_\eps$ of the Boltzmann equation (\ref{BoltzEqM}) is a small
perturbation of the uniform Maxwellian equilibrium state $M$. Therefore, the linearization about $M$ of the Boltzmann collision integral 
plays an important role in these limits. Thus, we consider this linearized collision integral intertwined with the multiplication by $M$, and set 
$\cL_M\phi=-M^{-1}D\cC(M)\cdot(M\phi)$, or equivalently
\begin{equation}\label{LinCollInt}
\cL_M\phi(v):=\iint_{\bR^3\times\bS^2}(\phi(v)+\phi(v_*)-\phi(v')-\phi(v'_*))|(v-v_*)\cdot\om|M(v_*)\dd v_*\dd\om\,.
\end{equation}
Hilbert \citep{Hilbert12} proved that $\cL_M$ is an unbounded, Fredholm, self-adjoint nonnegative operator on $L^2(\bR^3;M\dd v)$\footnote{The
notation $L^p(\bR^N;f\dd v)$ (where $p\ge 1$ and $f>0$ is a measurable function defined a.e. on $\Om$) designates the set of measurable 
functions $\phi$ defined a.e. on $\bR^N$ that satisfy
$$
\int_\Om|\phi(v)|^pf(v)\dd v<+\infty\,.
$$}, 
with domain $L^2(\bR^3;(1+|v|^2)M\dd v)$ and nullspace $\hbox{Ker}\cL_M=\hbox{Span}\{1,v_1,v_2,v_3,|v|^2\}$. Hilbert's argument, written for 
the hard sphere case, was later extended by Grad \citep{Grad63}, who defined some appropriate class of collision kernels $b(v-v_*,\om)$ for 
which the linearized collision integral satisfies the Fredholm alternative. Grad's idea was that grazing collisions between neutral gas molecules 
are rare events that can be somehow neglected, at variance with the case of plasmas or ionized gases. Henceforth, we denote
$$
\la\phi\ra=\int_{\bR^3}\phi(v)M(v)\dd v\qquad\hbox{ for each }\phi\in L^1(\bR^3;M\dd v)\,.
$$
With this notation, the Fredholm alternative for the integral equation $\cL_Mf=S$ with unknown $f$ and source term $S\in L^2(\bR^3;M\dd v)$ 
can be stated as follows:

a) either $\la S\ra=\la Sv_1\ra=\la Sv_2\ra=\la Sv_3\ra=\la S|v|^2\ra=0$, in which case the integral equation has a unique solution $f$ satisfying
$$
f\in L^2(\bR^3;(1+|v|^2)M\dd v)\quad\hbox{ and }\quad\la f\ra=\la fv_1\ra=\la fv_2\ra=\la fv_3\ra=\la f|v|^2\ra=0\,,
$$
henceforth denoted $f=\cL_M^{-1}S$, or

b) there exists $\phi\in\hbox{Ker}\cL_M$ such that $\la S\phi\ra\not=0$, in which case the integral equation $\cL_Mf=S$ does not have any
solution in $L^2(\bR^3;(1+|v|^2)M\dd v)$. 

%%%%%%%%%%%%%%%%%%%%%%%%%%%%%%%%%%%%%%%%%%%%%%%%%%%%%%%%%%%%%%%%%%%%%%%%%
\subsection{The moment method for the Navier-Stokes limit: formal argument}
%%%%%%%%%%%%%%%%%%%%%%%%%%%%%%%%%%%%%%%%%%%%%%%%%%%%%%%%%%%%%%%%%%%%%%%%%

Define $g_\eps$ by the formula $F_\eps=M(1+\de_\eps g_\eps)$. If $F_\eps$ satisfies (\ref{BoltzEqM}) with $\Kn=\Sh=\de_\eps=\eps$, the 
relative fluctuation of distribution function $g_\eps$ satisfies
\begin{equation}
\label{FluctBoltz}
\eps\d_tg_\eps+v\cdot\grad_xg_\eps+\frac1\eps \cL_Mg_\eps=\cQ_M(g_\eps,g_\eps)\,,
\end{equation}
where $\cQ_M$ is the symmetric bilinear operator defined by $\cQ_M(\phi,\phi)=M^{-1}\cC(M\phi)$. Multiplying each side of (\ref{FluctBoltz}) by
$\eps$ and letting $\eps\to 0$ shows that, if $g_{\eps_n}\to g$ for some subsequence $\eps_n\to 0$, the limiting fluctuation $g$ is a ``local
Maxwellian state'', i.e. is of the form
\begin{equation}
\label{InfinitMaxw}
g(t,x,v)=\rho(t,x)+u(t,x)\cdot v+\th(t,x)\tfrac12(|v|^2-3)\,.
\end{equation}

Multiplying each side of (\ref{FluctBoltz}) by $M$ and $vM$, and integrating in $v\in\bR^3$ shows that
$$
\eps\d_t\la g_\eps\ra+\Div_x\la vg_\eps\ra=0\,,\qquad\eps\d_t\la v g_\eps\ra+\Div_x\la v\otimes vg_\eps\ra=0
$$
in view of (\ref{CollCons}). Passing to the limit as $\eps_n\to 0$, and taking into account the local Maxwellian form (\ref{InfinitMaxw}) of $g$ leads
to
\begin{equation}
\label{IncomprBouss}
\Div_xu=\Div_x\la vg\ra=0\,,\qquad\grad_x(\rho+\th)=\Div_x\la v\otimes vg\ra=0\,.
\end{equation}
The first equality is the solenoidal condition for the velocity field $u$, while the second implies that $\rho+\th=0$, assuming that $\rho,\th\to 0$
as $|x|\to+\infty$.

Next we multiply each side of (\ref{FluctBoltz}) by $\frac1\eps vM$ and integrate in $v\in\bR^3$ to obtain
$$
\d_t\la v g_\eps\ra+\Div_x\frac1\eps\la Ag_\eps\ra=-\grad_x\frac1\eps \La\frac{|v|^2}{3} g_\eps\Ra\,,
$$
where $A(v):=v\otimes v-\frac13 |v|^2I$. One has $\la A_{kl}\ra=\la A_{kl}v_1\ra=\la A_{kl}v_2\ra=\la A_{kl}v_3\ra=\la A_{kl}|v|^2\ra=0$ for
each $k,l=1,2,3$, so that $\hat A_{kl}:=\cL_M^{-1}A_{kl}\in L^2(\bR^3;(1+|v|^2)M\dd v)$ is well-defined. Since $\cL_M$ is self-adjoint on 
$L^2(\bR^3;M\dd v)$, one has
\begin{equation}
\label{Flux}
\begin{array}{rl}
\frac1\eps\la Ag_\eps\ra=\frac1\eps\la(\cL_M\hat A)g_\eps\ra=\La\hat A\frac1\eps\cL_Mg_\eps\Ra\!\!\!\!\!\!
&=\la\hat A(\cQ_M(g_\eps,g_\eps)\ra-\la\hat A(\eps\d_tg_\eps+v\cdot\grad_xg_\eps)\ra
\\
&\to\la\hat A(\cQ_M(g,g)\ra-\la\hat Av\cdot\grad_xg\ra\,.
\end{array}
\end{equation}
By (\ref{InfinitMaxw}) and the solenoidal condition in (\ref{IncomprBouss}), the second term takes the form
\begin{equation}
\label{ViscTerm}
\la\hat Av\cdot\grad_xg\ra=\la\hat A\otimes v\otimes v\ra:\grad_xu=\nu(\grad_xu+(\grad_xu)^T)\,.
\end{equation}
Indeed $\la\hat A_{ij}A_{kl}\ra=\nu(\de_{ik}\de_{jl}+\de_{il}\de_{jk}-\tfrac23 \de_{ij}\de_{kl})$, which can be recast as $\nu=\frac1{10} \la\hat A:A\ra$
since $A(Rv)=RA(v)R^T$ for each $v\in\bR^3$ and each $R\in O_3(\bR)$. This formula for $\nu$ is equivalent to the first relation in (\ref{nu-ka}).

As for the first term, since $g\in\hbox{Ker}\cL_M$ according to (\ref{InfinitMaxw}), one has $\cQ_M(g,g)=\frac12\cL_M(g^2)$ (see \citep{BGL1}, fla.
(60) on p. 338.) Hence
\begin{equation}
\label{NLTerm}
\la\hat A(\cQ_M(g,g)\ra=\tfrac12 \la\hat A\cL_M(g^2)\ra=\tfrac12 \la(\cL_M\hat A)g^2\ra=\tfrac12 \la Ag^2\ra
	=\tfrac12 \la A\otimes v\otimes v\ra:u\otimes u=u\otimes u-\tfrac13 |u|^2I\,,
\end{equation}
in view of the elementary identity $\la A_{ij}A_{kl}\ra=(\de_{ik}\de_{jl}+\de_{il}\de_{jk}-\tfrac23 \de_{ij}\de_{kl})$.

Let $\xi\equiv\xi(x)\in C^\infty_c(\bR^3)$ be a divergence-free test vector field. Substituting (\ref{ViscTerm}) and (\ref{NLTerm}) in (\ref{Flux}) shows 
that
$$
\begin{array}{rl}
\displaystyle 0\!\!\!\!\!\!&=\displaystyle\d_t\int\xi\cdot\la v g_\eps\ra\dd x-\int\grad\xi:\frac1\eps\la Ag_\eps\ra\dd x
\\	\\
&\displaystyle \to\d_t\int\xi\cdot u\dd x-\int\grad\xi:(u\otimes u-\tfrac13 |u|^2I)\dd x+\nu\int\grad\xi:(\grad_xu+(\grad_x u)^T)\dd x\,.
\end{array}
$$
Since $\xi$ is divergence-free
$$
\int\grad\xi:\tfrac13 |u|^2I\dd x=\tfrac13 \int|u|^2\Div_x\xi\dd x=0\,,
$$
while
$$
\int\grad\xi:(\grad_x u)^T\dd x=-\int\grad(\Div\xi)\cdot u\dd x=0\,.
$$
Therefore
\begin{equation}\label{WeakFormNS}
\d_t\int\xi\cdot u\dd x-\int\grad\xi:u\otimes u\dd x+\nu\int\grad\xi:\grad_xu\dd x=0
\end{equation}
for each divergence-free test vector field $\xi\equiv\xi(x)\in C^\infty_c(\bR^3)$. Now, if $T\in\cD'(\bR^3)$ is a vector-valued distribution
satisfying $\la T,\xi\ra=0$ for each divergence-free test vector field $\xi\equiv\xi(x)\in C^\infty_c(\bR^3)$, there exists a scalar distribution 
$\pi\in\cD'(\bR^3)$ such that $T=\grad\pi$. In other words, (\ref{WeakFormNS}) means precisely that $u$ is a weak solution of the motion
equation in the Navier-Stokes system.

%%%%%%%%%%%%%%%%%%%%%%%%%%%%%%%%%%%%%%%%%%%%%%%%%%%%%%%%%%%%%%%%%%%%%%%%%
\subsection{Compactness tools}
%%%%%%%%%%%%%%%%%%%%%%%%%%%%%%%%%%%%%%%%%%%%%%%%%%%%%%%%%%%%%%%%%%%%%%%%%

An important ingredient in the proof of all fluid dynamic limits of the Boltzmann equation considered above is the fact that the relative
fluctuation of distribution function $g_\eps:=(F_\eps-M)/\de_\eps M$ converges in some sense, possibly after extracting some subsequence
$\eps_n\to 0$. The key argument is the following inequality resulting from (\ref{EntrIneq}):
$$
\int_{\bR^3}\la h(\de_\eps g_\eps)\ra(t)\dd x=H(F_\eps(t,\cdot,\cdot)|M)\le H(F_\eps^{in}|M)=O(\de_\eps^2)\qquad\hbox{ for each }t\ge 0\,,
$$
for the initial data considered in the four theorems stated in the previous section, where $h(z):=(1+z)\ln(1+z)-z$. Since $h(z)\simeq z^2/2$ as 
$z\to 0$, this control is as good as a bound in $L^\infty(\bR_+;L^2(\bR^3\times\bR^3;M\dd x\dd v))$ for the values of $g_\eps$ not exceeding 
$O(1/\de_\eps)$. Thus  $(1+|v|^2)g_\eps$ is relatively compact in weak-$L^1([0,T]\times[-R,R]^3\times\bR^3; M\dd t\dd x\dd v)$ for each $R,T>0$, 
and  all its limit points as $\eps\to 0$ belong to $L^\infty(\bR_+;L^2(\bR^3\times\bR^3;M\dd x\dd v))$. In the case of the acoustic or Stokes-Fourier 
limit, the uniqueness of the solution of the limiting fluid equations implies that the whole family $(1+|v|^2)g_\eps$ converges weakly.

Since the leading order term in (\ref{FluctBoltz}) is $\frac1\eps \cL_Mg_\eps$ and $\cL_M$ is a linear operator, the weak compactness of the
family $(1+|v|^2)g_\eps$ is enough to conclude that any limit point $g$ of that family as $\eps\to 0$ must satisfy $\cL_Mg=0$, and therefore
is an infinitesimal Maxwellian, i.e. is of the form (\ref{InfinitMaxw}).

In addition, for the Navier-Stokes-Fourier limit, the compactness of the family $g_\eps$ in $L^1_{loc}$ for the strong topology (implying the a.e.
pointwise convergence of a subsequence) is needed to pass to the limit in nonlinear terms. We use repeatedly some compactness results for
moments of the distribution function in the velocity variable $v$ based on bounds on the streaming operator --- see \citep{GPS,GLPS}. These 
compactness results are referred to as compactness by \textit{velocity averaging}. A typical example of a velocity averaging theorem used in 
the Navier-Stokes limit of the Boltzmann equation is as follows. We state it in the steady case for the sake of simplicity.

\smallskip
\noindent
\textbf{Theorem.} (Golse, Saint-Raymond \citep{GSR2002})
Let $f_n\equiv f_n(x,v)$ be a bounded sequence in $L^1(\bR^N_x\times\bR^N_v)$ such that the sequence $v\cdot\grad_xf_n$ is bounded in 
$L^1(\bR^N_x\times\bR^N_v)$, while $f_n$ itself is bounded in $L^1(\bR^N_x;L^p(\bR^N_v))$ for some $p>1$. Then

\smallskip
\noindent
a) $f_n$ is weakly relatively compact in $L^1_{loc}(\bR^N_x\times\bR^N_v)$; and

\noindent
b) for each $\phi\in C_c(\bR^N)$, the sequence of velocity averages
$$
\int_{\bR^N}f_n(x,v)\phi(v)\dd v
$$
is strongly relatively compact in $L^1_{loc}(\bR^N)$.

%%%%%%%%%%%%%%%%%%%%%%%%%%%%%%%%%%%%%%%%%%%%%%%%%%%%%%%%%%%%%%%%%%%%%%%%%
\subsection{The conservation laws}
%%%%%%%%%%%%%%%%%%%%%%%%%%%%%%%%%%%%%%%%%%%%%%%%%%%%%%%%%%%%%%%%%%%%%%%%%

The formal argument presented above in the case of the Navier-Stokes limit shows the importance of the local conservation laws of mass, 
momentum and energy in the derivation of fluid dynamic models from the Boltzmann equation. Unfortunately, renormalized solutions of the 
Boltzmann equation are not known to satisfy the local conservation laws of momentum and energy in (\ref{ConsLaws}). They satisfy instead
the approximate conservation laws
\begin{equation}\label{ApproxConsLaws}
\Sh\,\d_t\int_{\bR^3}\Gamma\left(\frac{F_\eps}{M}\right)\left(\begin{array}{c}v \\ \tfrac12 |v|^2\end{array}\right)M\dd v
+
\Div_x\int_{\bR^3}\Gamma\left(\frac{F_\eps}{M}\right)\left(\begin{array}{c} v\otimes v\\ \tfrac12 |v|^2v\end{array}\right)M\dd v
=
\frac1{\eps} \int_{\bR^3}\Gamma'\left(\frac{F_\eps}{M}\right)\cC(F_\eps)\left(\begin{array}{c}v \\ \tfrac12 |v|^2\end{array}\right)\dd v\,.
\end{equation}
Therefore, one must show that the conservation defects
$$
\frac1{\eps\de_\eps\Sh} \int_{\bR^3}\Gamma'\left(\frac{F_\eps}{M}\right)\cC(F_\eps)\left(\begin{array}{c}v \\ |v|^2\end{array}\right)\dd v\to 0
$$
in $L^1_{loc}(\bR_+\times\bR^3)$ as $\eps\to 0$, and identify the limits as $\eps\to 0$ of the terms
$$
\frac1{\de_\eps}\int_{\bR^3}\left(\Gamma\left(\frac{F_\eps}{M}\right)-\Gamma(1)\right)\left(\begin{array}{c}v \\ \tfrac12 |v|^2\end{array}\right)M\dd v
\quad\hbox{ and }\quad
\frac1{\de_\eps\Sh}\int_{\bR^3}\left(\Gamma\left(\frac{F_\eps}{M}\right)-\Gamma(1)\right)
	\left(\begin{array}{c}v\otimes v-\tfrac13 |v|^2I\\ \\ (|v|^2-5)v\end{array}\right)M\dd v\,.
$$
This raises an important question regarding the tail of the distribution functions $F_\eps$ as $|v|\to+\infty$. That the family $(1+|v|^2)g_\eps$ 
is relatively compact in weak-$L^1([0,T]\times[-R,R]^3\times\bR^3; M\dd t\dd x\dd v)$ for each $R,T>0$ is in general not enough --- for instance,
in the acoustic limit, $\Sh=1$ and one needs to identify the limit of the energy flux
$$
\frac1{\de_\eps}\int_{\bR^3}\left(\Gamma\left(\frac{F_\eps}{M}\right)-\Gamma(1)\right)\tfrac12|v|^2vM\dd v\,,
$$ 
which is a 3rd order moment in $v$ of $\frac1{\de_\eps} (\Gamma(1+\de_\eps g_\eps)-\Gamma(1))\sim \Gamma'(1)g_\eps$. Controlling the
high speed tail of (fluctuations of) the distribution function is an essential step in the derivation of fluid dynamic limits of the Boltzmann
equation, and involves rather technical estimates based on the entropy and entropy production estimates (\ref{EntrIneq}) together with the
dispersion effects of the streaming operator $\Sh\,\d_t+v\cdot\grad_x$ (see \citep{BGL2,GSR2004,GSR2009}).

%%%%%%%%%%%%%%%%%%%%%%%%%%%%%%%%%%%%%%%%%%%%%%%%%%%%%%%%%%%%%%%%%%%%%%%%%
\subsection{The relative entropy method}
%%%%%%%%%%%%%%%%%%%%%%%%%%%%%%%%%%%%%%%%%%%%%%%%%%%%%%%%%%%%%%%%%%%%%%%%%

In inviscid hydrodynamic limits, i.e. the compressible or incompressible Euler limits, entropy production does not balance streaming. Therefore
the velocity averaging method fails for such limits. The idea is to use the regularity of the solution of the target equations, together with relaxation 
towards local equilibrium in order to obtain some compactness on fluctuations of the distribution function.

\smallskip
Pick for instance $u$, a smooth solution of the target equations --- e.g. the incompressible Euler equations --- and study the evolution of the quantity
$$
Z_\eps(t)=\frac1{\de_\eps^2}H(F_\eps|\cM_{(1,\de_\eps u(t,x),1)})
$$
where $F_\eps$ is a renormalized solution of (\ref{BoltzEqM}) with $\Sh=\Ma=\de_\eps=\eps^\a$ and $0<\a<1$. This is the leading order of the 
relative entropy of the Boltzmann solution with respect to the local Maxwellian state defined by $u$, in the incompressible Euler scaling. At the 
formal level, it is found that
$$
\frac{\dd Z_\eps}{\dd t}(t)=-\frac1{\de_\eps^2}\int_{\bR^3}\grad_xu:\int_{\bR^3}(v-\de_\eps u)^{\otimes 2}F_\eps\dd v\dd x
+\frac1{\de_\eps}\int_{\bR^3}\grad_xp\cdot\int_{\bR^3}(v-\de_\eps u)F_\eps\dd v\dd x\,.
$$
The second term in the r.h.s. of the equality above vanishes with $\eps$ since
$$
\frac1{\de_\eps}\int_{\bR^3}vF_\eps(t,x,v)\dd v\to\hbox{ divergence-free vector field.}
$$
The key idea is to estimate the first term in the r.h.s. as follows
$$
\frac1{\de_\eps^2}\iint_{\bR^3\times\bR^3}\left|\grad_xu:(v-\de_\eps u)^{\otimes 2}F_\eps\right|\dd v\dd x\le CZ_\eps(t)+o(1)
$$
where $C=O(\|\grad_xu\|_{L^\infty})$. Then, one concludes with Gronwall's inequality.

\smallskip
The relative entropy method stems from an idea of H.T. Yau (for Ginzburg-Landau lattice models, see \citep{HTYau93}); it was later adapted to 
the Boltzmann equation by the author \citep{BouGolPul2000} and Lions-Masmoudi \citep{LionsMasmoudi2}. It is especially designed to handle
sequences of \textit{weak solutions of the Boltzmann equation} converging to a \textit{classical solution of the fluid equation}.

%%%%%%%%%%%%%%%%%%%%%%%%%%%%%%%%%%%%%%%%%%%%%%%%%%%%%%%%%%%%%%%%%%%%%%%%%
\section{Conclusions}
%%%%%%%%%%%%%%%%%%%%%%%%%%%%%%%%%%%%%%%%%%%%%%%%%%%%%%%%%%%%%%%%%%%%%%%%

The DiPerna-Lions theory of renormalized solutions of the Boltzmann equation allows one to obtain derivations of fluid dynamic regimes from 
the kinetic theory of gases without unphysical assumptions on the size or regularity of the data. Following the program outlined in \citep{BGL2},  
these derivations are based on

\noindent
a) relative entropy and entropy production estimates, together with

\noindent
b) functional analytic methods in Lebesgue ($L^p$) spaces.

%\smallskip
At present, the program in \citep{BGL2} leaves aside the compressible Euler limit of the Boltzmann equation, or the asymptotic regime leading to 
the compressible Navier-Stokes equations. Little progress has been made on these issues since the work of Nishida \citep{Nishida78} and
Caflisch \citep{Caflisch80}.

%\smallskip
The problem of deriving fluid dynamic limits from the Boltzmann equation in the \textit{steady} regime is also of considerable importance for
practical applications. Formal results are of course well understood with the classical Hilbert or Chapman-Enskog expansion techniques --- see
the book of Sone \citep{SoneBook1}. Unfortunately, the theory of the steady Boltzmann equation with large data is not as mature as its counterpart 
for the Cauchy problem, in spite of interesting contributions by Arkeryd and Nouri \citep{ArkerydNouri}, and there is no analogue of the DiPerna-Lions 
theory for the steady case yet.

%\smallskip
But even for evolution problems in regimes that are weakly nonlinear at the kinetic level, the relative entropy is not the solution to all difficulties. 
In several asymptotic regimes of the Boltzmann equation, the leading order and next to leading order fluctuations of the distribution function may
interact to produce highly nontrivial macroscopic effects in the fluid dynamic limit. Examples of such asymptotic regimes are

%\smallskip
\noindent
a) ghost effects, introduced by Sone, Aoki, Takata, Sugimoto and Bobylev in \citep{SoneAokiGhost}, reported in Sone's Harold Grad Lecture 
\citep{SoneGrad} and \citep{SoneBook1,SoneBook2}, and systematically studied by Sone, Aoki and the Kyoto school,

\noindent
b) Navier-Stokes limits recovering viscous heating terms, due to Bobylev \citep{Bobylev95} and Bardos-Levermore-Ukai-Yang
 \citep{BardLvrmrUkaiYang2009} --- see also the discussion in \citep{SoneBardFGSugi2000}, and

\noindent
c) hydrodynamic limits for thin layers of fluid --- see \citep{Golse2010}.

\smallskip
\noindent
\textbf{Acknowledgements.} The author thanks Profs. Aoki, Levermore and Sone for their generous scientific advice during the preparation of this 
paper.

%%%%%%%%%%%%%%%%%%%%%%%%%%%%%%%%%%%%%%%%%%%%%%%%%%%%%%%%%%%%%%%%%%%%%%%%


\begin{thebibliography}{99}

\bibitem[Arkeryd et al. (2002)]{ArkerydNouri}
L. Arkeryd, and A. Nouri, \emph{Ann. Sc. Norm. Super. Pisa Cl. Sci. (5)}  {\bf 1}  (2002), pp. 359--385.

\bibitem[Arkeryd et al. (2006)]{ArkNouCouette}
L. Arkeryd, and A. Nouri, \emph{J. Stat. Phys.}  {\bf 124}  (2006), pp. 401--443.

\bibitem[Bardos et al. (1991)]{BGL1}
C. Bardos, F. Golse, and C.D. Levermore, \emph{J. Stat. Phys.} {\bf 63} (1991), pp. 323--344.

\bibitem[Bardos et al. (1993)]{BGL2}
C. Bardos, F. Golse, and C.D. Levermore, \emph{Comm. Pure \& Appl. Math} {\bf 46} (1993), pp. 667--753.

\bibitem[Bardos et al. (2000)]{BGL3}
C. Bardos, F. Golse, and C.D. Levermore, \emph{Arch. Ration. Mech. Anal.} {\bf 153} (2000), pp. 177--204.

\bibitem[Bardos et al. (2009)]{BardLvrmrUkaiYang2009}
C. Bardos, C.D. Levermore, S. Ukai, and T. Yang, \emph{Bull. Inst. Math. Acad. Sin. (N.S.)} \textbf{3} (2008), pp. 1--49.

\bibitem[Bardos et al. (1991)]{BardosUkai91}
C. Bardos, and S. Ukai, \emph{Math. Models and Methods in the Appl. Sci.} {\bf 1} (1991), pp. 235--257.

\bibitem[Bobylev (1995)]{Bobylev95}
A. Bobylev, \emph{J. Statist. Phys.} {\bf 80} (1995), pp. 1063--1083. 

\bibitem[Bouchut et al. (2000)]{BouGolPul2000}
F. Bouchut, F. Golse, and M. Pulvirenti, \emph{Kinetic Equations and Asymptotic Theory}, edited by L. Desvillettes \&
B. Perthame, Editions scientifiques et m\'edicales Elsevier, Paris, 2000.

\bibitem[Caflisch (1980)]{Caflisch80}
R. Caflisch, \emph{Comm. on Pure and Appl. Math.} {\bf 33} (1980), pp. 651--666.

\bibitem[DeMasi et al. (1990)]{dMEL}
A. DeMasi, R. Esposito, and J. Lebowitz, \emph{Comm. on Pure and Appl. Math.} {\bf 42} (1990), pp. 1189--1214.

\bibitem[DiPerna et al. (1990)]{diPernaLions90}
R. DiPerna, and P.-L. Lions, \emph{Ann. of Math.} {\bf 130} (1990), pp. 321--366.

\bibitem[Golse (2010)]{Golse2010}
F. Golse, \emph{From the Boltzmann equation to fluid dynamics in thin layers}, preprint (2010).

\bibitem[Golse et al. (2002)]{GolseLvrmr2002}
F. Golse, and C.D. Levermore, \emph{Comm. on Pure and Appl. Math.} {\bf 55}  (2002), pp. 336--393.

\bibitem[Golse et al. (1988)]{GLPS}
F. Golse, P.-L. Lions, B. Perthame, R. Sentis, \emph{J. Func. Anal.} {\bf 76} (1988), pp. 110--125.

\bibitem[Golse et al. (1985)]{GPS}
F. Golse, B. Perthame, R. Sentis, \emph{C.R. Acad. Sci. S\'er. I}  {\bf 301} (1985), pp. 341--344.

\bibitem[Golse et al. (2002)]{GSR2002}
F. Golse, and L. Saint-Raymond, \emph{C. R. Acad. Sci. S\'er. I Math.} {\bf 334} (2002), pp. 557--562.

\bibitem[Golse et al. (2004)]{GSR2004}
F. Golse, and L. Saint-Raymond, \emph{Invent. Math.}  {\bf 155}  (2004),  pp. 81--161.

\bibitem[Golse et al. (2009)]{GSR2009}
F. Golse, and L. Saint-Raymond, \emph{J. Math. Pures et Appl.} {\bf 91} (2009), 508--552.

\bibitem[Grad (1962)]{Grad63}
H. Grad, in Rarefied Gas Dynamics, edited by J. A. Laurmann (Academic Press, New York, 1963), Vol. 1, pp. 26-59.

\bibitem[Hilbert (1900)]{Hilbert00}
D. Hilbert, in Internat. Congress of Math., Paris 1900, transl. and repr. in \emph{Bull. Amer. Math. Soc.} {\bf 37} (2000), pp. 407-436.

\bibitem[Hilbert (1912)]{Hilbert12}
D. Hilbert, \emph{Math. Ann.} {\bf 72} (1912), pp. 562--577.

\bibitem[Jing et al. (2009)]{JingLvrmrMasmoudi}
N. Jiang, C.D. Levermore, and N. Masmoudi, preprint arxiv 0903.5086.

\bibitem[Leray (1934)]{Leray1934}
J. Leray, \emph{Acta Math.} {\bf 63} (1934), pp. 193--248.

\bibitem[Levermore et al. (2010)]{LvrmrMasmoudi2010}
C.D. Levermore, and N. Masmoudi, \emph{Archive Rat. Mech. \& Anal.} {\bf 196} (2010), pp. 753--809.

\bibitem[Lions (1993)]{LionsKyoto93}
P.-L. Lions, \emph{J. Math. Kyoto Univ.} {\bf 34} (1994), 429--461.

\bibitem[Lions (1994)]{Lions94}
P.-L. Lions, \emph{Comm. in Partial Differential Equations} {\bf 19} (1994), pp. 335--367.

\bibitem[Lions et al. (2001)]{LionsMasmoudi1}
P.-L. Lions, and N. Masmoudi, \emph{Archive Rat. Mech. \& Anal.} {\bf 158} (2001), pp. 173--193.

\bibitem[Lions et al. (2001)]{LionsMasmoudi2}
P.-L. Lions, and N. Masmoudi, \emph{Archive Rat. Mech. \& Anal.} {\bf 158} (2001), pp. 195--211.

\bibitem[Maxwell (1866)]{Maxwell66}
J. Clerk Maxwell, \emph{Philosophical Transactions} {\bf 157} (1866).

\bibitem[Nishida (1978)]{Nishida78}
T. Nishida, \emph{Comm. Math. Phys.}  {\bf 61}  (1978), pp. 119--148.

\bibitem[Saint-Raymond (2003)]{SR2003}
L. Saint-Raymond, \emph{Arch. Ration. Mech. Anal.} {\bf 166} (2003), pp. 47--80.

\bibitem[Sideris (1987)]{Sideris87}
T. Sideris, \emph{Commun. Math. Phys.} {\bf 101} (1985), pp. 475--485.

\bibitem[Sone et al. (1996)]{SoneAokiGhost}
Y. Sone, K. Aoki, S. Takata, H. Sugimoto, and A. V. Bobylev, \emph{Phys. Fluids} {\bf 8} (1996), pp. 628-638.

\bibitem[Sone (1997)]{SoneGrad}
Y. Sone, in Rarefied Gas Dynamics, edited by C. Shen (Peking University Press, Beijing, 1997), pp. 3--24.

\bibitem[Sone (2002)]{SoneBook1}
Y. Sone, \emph{Kinetic Theory and Fluid Dynamics}, Birkh\"auser, Boston, 2002.

\bibitem[Sone (2007)]{SoneBook2}
Y. Sone, \emph{Molecular Gas Dynamics, Theory, Techniques, and Applications}, Birkh\"auser, Boston, 2007.

\bibitem[Sone et al. (2000)]{SoneBardFGSugi2000}
Y. Sone, C. Bardos, F. Golse, H. Sugimoto, \emph{Eur. J. Mech. B - Fluids} {\bf 19} (2000), pp. 325--360.

\bibitem[Yau (1993)]{HTYau93}
H.T. Yau, \emph{Lett. Math. Phys.}  {\bf 22}  (1991), pp. 63--80.

\end{thebibliography}
\end{document}